\documentclass{article}

\usepackage{amsmath}
\usepackage{amssymb}

\newtheorem{thm}{Theorem}
\newtheorem{conj}{Conjecture}
\newtheorem{lem}{Lemma}
\newtheorem{ques}{Question}

\title{Approximating reals by sums of two rationals}
\author{Tsz Ho Chan}

\begin{document}
\maketitle
\begin{abstract}
We generalize Dirichlet's diophantine approximation theorem to
approximating any real number $\alpha$ by a sum of two rational
numbers $\frac{a_1}{q_1} + \frac{a_2}{q_2}$ with denominators $1
\leq q_1, q_2 \leq N$. This turns out to be related to the
congruence equation problem $x y \equiv c \pmod q$ with $1 \leq x, y
\leq q^{1/2 + \epsilon}$.
\end{abstract}

\section{Introduction}

Dirichlet's theorem on rational approximation says
\begin{thm} \label{thm1}
For any real number $\alpha$ and real number $N \geq 1$, there exist
integers $1 \leq q \leq N$ and $a$ such that
\begin{equation} \label{dirichlet}
\Big|\alpha - \frac{a}{q}\Big| \leq \frac{1}{q N}.
\end{equation}
\end{thm}

While studying almost squares (see [\ref{C}]), the author
accidentally consider the question of approximating $\alpha$ by a
sum of two rational numbers:
\begin{ques} \label{q1}
Find a good upper bound for
\begin{equation} \label{d2}
\Big| \alpha - \frac{a_1}{q_1} - \frac{a_2}{q_2} \Big|
\end{equation}
with integers $a_1, a_2$ and $1 \leq q_1, q_2 \leq N$.
\end{ques}
(This turns out to be unfruitful towards the study of almost
squares.) One can continue further with approximating $\alpha$ by a
sum of $n > 2$ rational numbers. We shall study this in another
paper. These seem to be some new questions in diophantine
approximation.

However, if one combines the two fractions in (\ref{d2}), it becomes
$$\Big| \alpha - \frac{b}{q_1 q_2} \Big|$$
with integers $b$ and $1 \leq q_1, q_2 \leq N$. This looks like the
left hand side of (\ref{dirichlet}) except that we require the
denominator $q$ to be of a special form, namely $q = q_1 q_2$ with
$1 \leq q_1, q_2 \leq N$. In this light, our Question \ref{q1} is
not so new after all. People have studied diophantine approximation
where the denominator $q$ is of some special form. For example,
\begin{enumerate}
\item When $q$ is a perfect square, see Zaharescu [\ref{Z}] with some
history of the problem.
\item When $q$ is squarefree, see Harman [\ref{Ha1}], Balog \& Perelli [\ref{BP}],
Heath-Brown [\ref{He1}].
\item When $q$ is prime, see Heath-Brown \& Jia [\ref{HJ}] with some
history of the problem.
\item When $q$ is B-free, see Alkan, Harman and Zaharescu [\ref{A}].
\end{enumerate}
Techniques from exponential sum, character sum, sieve method, and
geometry of numbers were used in the above list of works. In this
paper, we shall use exponential sum and character sum methods to
study Question \ref{q1}. It would be interesting to see if other
methods can be applied. One distinct feature of our results (see
next section) is that the upper bounds of (\ref{d2}) depend on
single rational approximations $\frac{a}{q}$ of the real number
$\alpha$ given by Dirichlet's Theorem. Alternatively, we try to see
how approximation by a sum of two rationals compares with single
rational approximation.

The starting point of the argument is
\begin{equation} \label{d3}
\Big| \alpha - \frac{a_1}{q_1} - \frac{a_2}{q_2} \Big| \leq \Big|
\alpha - \frac{a}{q} \Big| + \Big| \frac{a}{q} - \frac{a_1}{q_1} -
\frac{a_2}{q_2} \Big|
\end{equation}
by triangle inequality. The first term on the right hand side of
(\ref{d3}) is small by Dirichlet's Theorem. Thus, it remains to
obtain a good upper bound for the second term on the right hand side
of (\ref{d3}) (i.e. we can restrict Question \ref{q1} to rational
$\alpha$). By combining denominators and letting $b = a_1 q_2 + a_2
q_1$,
$$\Big| \frac{a}{q} - \frac{a_1}{q_1} - \frac{a_2}{q_2} \Big| =
\frac{|a q_1 q_2 - b q|}{q q_1 q_2}.$$ Our goal is trying to make
the numerator as small as possible, say $a q_1 q_2 - b q = r$ where
$r > 0$ is small. Transforming this into a congruence equation
$\pmod q$, we have
\begin{equation} \label{d4}
q_1 q_2 \equiv r \overline{a} \pmod q
\end{equation}
with $1 \leq q_1, q_2 \leq N$ where $a \overline{a} \equiv 1 \pmod
q$. This gives some indication why Question \ref{q1} is related to
the congruence equation problem stated in the abstract. Ideally, we
want to solve (\ref{d4}) with $r = 1$. This seems too hard. So, we
take advantage of allowing $r$ to run over a short interval which
gives Theorems \ref{thm6}, \ref{thm7} and \ref{thm5} in the next
section. We will also prove an almost all result, namely Theorem
\ref{thm4}, towards Conjecture \ref{conj2}.

Throughout the paper, $\epsilon$ denotes a small positive number.
Both $f(x) = O(g(x))$ and $f(x) \ll g(x)$ mean that $|f(x)| \leq C
g(x)$ for some constant $C > 0$. Moreover $f(x) = O_\lambda(g(x))$
and $f(x) \ll_\lambda g(x)$ mean that the implicit constant $C =
C_\lambda$ may depend on the parameter $\lambda$. Also $\phi(n)$ is
Euler's phi function and $d(n)$ is the number of divisors of $n$.
Finally $|\mathcal{S}|$ stands for the cardinality of the set $S$.
\section{Some conjectures and results}
By imitating Dirichlet's theorem, one might conjecture that there
exist integers $a_1$, $a_2$, $1 \leq q_1, q_2 \leq N$ such that
$$\Big|\alpha - \frac{a_1}{q_1} - \frac{a_2}{q_2}\Big| \ll \frac{1}{q_1 q_2 N^2}$$
as the two fractions combine to give a single fraction with
denominator $q_1 q_2 \leq N^2$. However, this is very wrong as
illustrated by the following example:

Let $\alpha  = \frac{a}{p}$ for some prime number $p$ with $N < p
\leq 2N$ (guaranteed to exist by Bertrand's postulate) and integer
$a$ with $(a,p)=1$. Then
\begin{equation} \label{example}
\Big|\alpha - \frac{a_1}{q_1} - \frac{a_2}{q_2}\Big| =
\Big|\frac{a}{p} - \frac{a_1}{q_1} - \frac{a_2}{q_2}\Big| \geq
\frac{1}{p q_1 q_2} \gg \frac{1}{q_1 q_2 N}.
\end{equation}

So, the best upper bound one can hope for is
\begin{equation} \label{bound1}
\Big|\alpha - \frac{a_1}{q_1} - \frac{a_2}{q_2}\Big| \ll
\frac{1}{q_1 q_2 N}
\end{equation}
for some integers $a_1$, $a_2$, $1 \leq q_1, q_2 \leq N$. This
follows directly from Theorem \ref{thm1} by simply choosing
$\frac{a_2}{q_2} = \frac{0}{1}$. On the other hand, one can easily
get the bound
\begin{equation} \label{bound2}
\Big|\alpha - \frac{a_1}{q_1} - \frac{a_2}{q_2}\Big| \ll
\frac{1}{N^2}.
\end{equation}
For example, fix $q_1$ and $q_2$ to be two distinct primes in the
interval $[N/4, N]$ (without loss of generality, we may assume $N
\geq 12$), and consider the fractions $\frac{k}{q_1 q_2}$ with $(k,
q_1 q_2) = 1$.  Since one of $k$, $k+1$ or $k+2$ is not divisible by
neither $q_1$ nor $q_2$, the distance between successive fractions
is $\leq \frac{3}{q_1 q_2} \ll \frac{1}{N^2}$. Interpolating between
(\ref{bound1}) and (\ref{bound2}), we make the following
\begin{conj} \label{conj0}
Let  $0 \leq \beta \leq 1$. For any real number $\alpha$ and real
number $N \geq 1$, there exist integers $1 \leq q_1, q_2 \leq N$ and
$a_1$ ,$a_2$ such that
$$\Big|\alpha - \frac{a_1}{q_1} - \frac{a_2}{q_2} \Big| \leq
\frac{1}{(q_1 q_2)^\beta N^{2 - \beta}}.$$
\end{conj}

From the above discussion, Conjecture \ref{conj0} is true when
$\beta = 0$ or $1$. We leave the cases $0 < \beta < 1$ to the
readers as a challenging open problem. In another direction, as
shown by the example in (\ref{example}), the approximation of a real
number by a sum of two rationals depends on the rational
approximation of $\alpha$ by a single rational number. Thus, we come
up with

\begin{conj} \label{conj1}
For any small $\epsilon > 0$ and any $N \geq 1$, suppose $\alpha$
has a rational approximation $|\alpha - \frac{a}{q}| \leq \frac{1}{q
N^2}$ for some integers $a$, $1 \leq q \leq N^2$ and $(a,q) = 1$.
Then
$$\Big|\alpha - \frac{a_1}{q_1} - \frac{a_2}{q_2}\Big| \ll_\epsilon
\frac{1}{q N^{2 - \epsilon}}$$ for some integers $a_1$, $a_2$, $1
\leq q_1, q_2 \leq N$. Note: We may restrict our attention to $q >
N$ in the above rational approximation of $\alpha$, for otherwise we
can just pick $\frac{a_1}{q_1} = \frac{a}{q}$ and $\frac{a_2}{q_2} =
\frac{0}{1}$.
\end{conj}

Roughly speaking, this means that if one can approximate a real
number well by a rational number, then one should be able to
approximate it by a sum of two rational numbers nearly as well. Note
that the example in (\ref{example}) shows that this conjecture is
best possible (apart from $\epsilon$). As indicated in the
Introduction, Conjecture \ref{conj1} is related to a conjecture on
congruence equation:
\begin{conj} \label{conj2}
Let $\epsilon$ be any small positive real number. For any positive
integer $q$ and integer $c$ with $(c,q) = 1$, the equation
$$x y \equiv c \pmod q$$
has solutions in $1 \leq x, y \ll_\epsilon q^{1/2 + \epsilon}$.
\end{conj}

We consider the following variation.
\begin{conj} \label{conj3}
Let $1/2 < \theta \leq 1$. There is a constant $C_\theta$ such that,
for any positive integer $q$ and integer $c$ with $(c,q) = 1$, the
equation
$$x y \equiv c \pmod q$$
has solutions in $ C_\theta N \leq x, y \leq 2 C_\theta N$ with
$(x,y) = 1$ for every $N \geq q^\theta$.
\end{conj}

Assume Conjecture \ref{conj3} for some $1/2 < \theta \leq 1$. Let
$|\alpha - \frac{a}{q}| \leq \frac{1}{q N^{1/\theta}}$ with
$(a,q)=1$ and $N < q \leq N^{1/\theta}$. Then, by Conjecture
\ref{conj3}, there are $C_\theta$, $q_1$, $q_2$ such that $C_\theta
N \leq q_1, q_2 \leq 2 C_\theta N$ and $q_1 q_2 \equiv \overline{a}
\pmod q$ (because $q^\theta \leq N$). Here $\overline{a}$ denotes
the multiplicative inverse of $a$ modulo $q$. In particular, we have
$a q_1 q_2 \equiv 1 \pmod q$. So, $a q_1 q_2 = k q + 1$ for some
integer $k$. This gives $a q_1 q_2 - k q = 1$ and
$$\Big|\frac{a}{q} - \frac{k}{q_1 q_2}\Big| = \frac{1}{q q_1 q_2}
\ll_\theta \frac{1}{q N^2}.$$ Since $(q_1, q_2) = 1$, the fraction
$\frac{k}{q_1 q_2} = \frac{a_1}{q_1} + \frac{a_2}{q_2}$ for some
integers $a_1$, $a_2$. Hence
$$\Big|\alpha - \frac{a_1}{q_1} - \frac{a_2}{q_2}\Big| \leq
\Big|\alpha - \frac{a}{q}\Big| + \Big|\frac{a}{q} - \frac{a_1}{q_1}
- \frac{a_2}{q_2} \Big| \ll_\theta \frac{1}{q N^{1/\theta}}.$$

Thus, if Conjecture \ref{conj3} is true with $\theta = \frac{1}{2} +
\epsilon$ for any small $\epsilon > 0$, then we have Conjecture
\ref{conj1}. From Igor Shparlinski (see [\ref{Sh}]), the author
learns that Conjecture \ref{conj3} is true for $\theta >
\frac{3}{4}$.

Recently, M.Z. Garaev and A.A. Karatsuba [\ref{gk}] proved that
\begin{thm}
Let $\Delta = \Delta(m) \rightarrow \infty$ as $m \rightarrow
\infty$. Then the set
$$\{xy \pmod m: 1 \leq x \leq m^{1/2}, S+1 \leq y \leq S + \Delta
m^{1/2} \sqrt{\frac{m}{\phi(m)}} \log m\}$$ contains $(1 +
O(\Delta^{-1})) m$ residue classes modulo $m$.
\end{thm}

One can think of this as an almost all result towards Conjecture
\ref{conj2}. With slight modification, one can get
\begin{thm} \label{thm3}
Let $\Delta = \Delta(m) \rightarrow \infty$ as $m \rightarrow
\infty$. For any small $\epsilon > 0$ and $m^{1/(2 - \epsilon)} \leq
N \leq m$, the set
$$\{xy \pmod m: \frac{N}{4} \leq x \leq N, S+1 \leq y \leq S + \Delta
m^{1/2} \sqrt{\frac{m}{\phi(m)}} \log m, (x,y) = 1\}$$ contains $(1
+ O_\epsilon(\Delta^{-1}) + O_\epsilon(\frac{\Delta}{\log m}
\sqrt{\frac{m}{\phi(m)}})) m$ residue classes modulo $m$.
\end{thm}

Using Theorem \ref{thm3}, we can prove an almost all result towards
Conjecture \ref{conj1}.
\begin{thm} \label{thm4}
Let $N$ be a positive integer. For any small $\epsilon > 0$,
Conjecture \ref{conj1} is true for all $0 \leq \alpha < 1$ except a
measure of $O_\epsilon(\frac{1}{\sqrt{\log N}})$.
\end{thm}

Instead of an almost all result, one may try to prove Conjecture
\ref{conj1} with a bigger uniform upper bound. Inspired by a recent
paper of Alkan, Harman and Zaharescu [\ref{A}], a variant of the
Erd\"{o}s-Turan inequality gives
\begin{thm} \label{thm6}
For any $\epsilon > 0$ and any $N \geq 1$, suppose $\alpha$ has a
rational approximation $|\alpha - \frac{a}{q}| \leq \frac{1}{q
N^{5/4}}$ for some integers $a$, $1 \leq q \leq N^{5/4}$ and $(a,q)
= 1$. Then
\begin{equation} \label{1thm6}
\Big|\alpha - \frac{a_1}{q_1} - \frac{a_2}{q_2}\Big| \ll_\epsilon
\frac{1}{q N^{5/4 - \epsilon}}
\end{equation}
for some integers $a_1$, $a_2$, and prime numbers $1 \leq q_1 < q_2
\leq N$.
\end{thm}
By a slight modification of the proof of Theorem \ref{thm6}, we have
\begin{thm} \label{thm7}
For any $\epsilon > 0$ and any $N \geq 1$, suppose $\alpha$ has a
rational approximation $|\alpha - \frac{a}{q}| \leq \frac{1}{q
N^{3/2}}$ for some integers $a$, $1 \leq q \leq N^{3/2}$ and $(a,q)
= 1$. Then
\begin{equation} \label{1thm7}
\Big|\alpha - \frac{a_1}{q_1} - \frac{a_2}{q_2}\Big| \ll_\epsilon
\frac{1}{q N^{3/2 - \epsilon}}
\end{equation}
for some integers $a_1$, $a_2$, $1 \leq q_1, q_2 \leq N$.
\end{thm}
Conditionally, we have
\begin{thm} \label{thm5}
Assume the Generalized Lindel\"{o}f Hypothesis. For any $\epsilon >
0$ and any $N \geq 1$, suppose $\alpha$ has a rational approximation
$|\alpha - \frac{a}{q}| \leq \frac{1}{q N^{4/3}}$ for some integers
$a$, $1 \leq q \leq N^{4/3}$ and $(a,q) = 1$. Then
\begin{equation} \label{2thm5}
\Big|\alpha - \frac{a_1}{q_1} - \frac{a_2}{q_2}\Big| \ll_\epsilon
\frac{1}{q N^{4/3 - \epsilon}}
\end{equation}
for some integers $a_1$, $a_2$, and prime numbers $1 \leq q_1 < q_2
\leq N$.
\end{thm}
This is proved by a character sum method. Of course, the goal is
trying to push the exponent of $N$ in (\ref{1thm6}), (\ref{1thm7})
and (\ref{2thm5}) to $2 - \epsilon$. Moreover, one may guess that
Conjecture \ref{conj1} is true even restricting $q_1$ and $q_2$ to
prime numbers.

The paper is organized as follow. We first prove the almost all
results, Theorems \ref{thm3} and \ref{thm4}, in section \ref{sec2}.
Then we prove Theorems \ref{thm6} and \ref{thm7} in sections
\ref{sec4} and \ref{sec5} respectively. Finally, we prove Theorem
\ref{thm5} in the last section.
\section{Almost all results: Theorem \ref{thm3} and \ref{thm4}} \label{sec2}

First, we assume Theorem \ref{thm3} and prove Theorem \ref{thm4}.

\bigskip

Proof of Theorem \ref{thm4}: Consider $1 \leq q \leq
N^{2-\epsilon}$, $0 \leq a \leq q$ with $(a,q) = 1$, and the
intervals $I_{a,q} = \{\alpha: |\alpha - \frac{a}{q}| \leq
\frac{1}{q N^{2-\epsilon}} \}$. By Theorem \ref{thm1}, $I_{a,q}$'s
cover the interval $[0,1]$. From the note in Conjecture \ref{conj1},
we may restrict our attention to $q
> N$.

We call $I_{a,q}$ good if there exist $\frac{N}{4} \leq q_1, q_2
\leq N$ with $(q_1, q_2) = 1$ such that $q_1 q_2 \equiv \overline{a}
\pmod q$. Otherwise, we call $I_{a,q}$ bad. By Theorem \ref{thm3},
for a fixed $q$, there are at most
\begin{equation} \label{bad}
\Bigl(O_\epsilon(\Delta^{-1}) + O_\epsilon \bigl(\frac{\Delta}{\log
q} \sqrt{\frac{q}{\phi(q)}}\bigr)\Bigr) q =
O\Bigl(\frac{q}{\sqrt{\log q}}
\bigl(\frac{q}{\phi(q)}\bigr)^{1/4}\Bigr)
\end{equation}
bad $I_{a,q}$'s by choosing $\Delta = \sqrt{\log q}
(\frac{\phi(q)}{q})^{1/4}$. For good $I_{a,q}$'s, we have $a q_1 q_2
- b q = 1$ for some integer $b$, which gives
$$\Big|\frac{a}{q} - \frac{a_1}{q_1} - \frac{a_2}{q_2}\Big| =
\Big|\frac{a}{q} - \frac{b}{q_1 q_2}\Big| = \frac{1}{q q_1 q_2} \ll
\frac{1}{q N^2}$$ for some integers $a_1$, $a_2$ as $(q_1, q_2) =
1$. Therefore, for $\alpha$ in good $I_{a,q}$, there exist
$\frac{N}{4} \leq q_1, q_2 \leq N$ and integers $a_1$, $a_2$ such
that
\begin{equation} \label{good}
\Big|\alpha - \frac{a_1}{q_1} - \frac{a_2}{q_2}\Big| \leq
\Big|\alpha - \frac{a}{q}\Big| + \Big|\frac{a}{q} - \frac{a_1}{q_1}
- \frac{a_2}{q_2}\Big| \ll_\epsilon \frac{1}{q N^{2-\epsilon}}
\end{equation}
Consequently by (\ref{bad}), (\ref{good}) holds for all but a
measure
$$\ll_\epsilon \sum_{N < q \leq N^{2-\epsilon}} \frac{q}{\sqrt{\log q}}
\bigl(\frac{q}{\phi(q)}\bigr)^{1/4} \times \frac{1}{q
N^{2-\epsilon}} \ll \frac{1}{\sqrt{\log N}}$$ of $\alpha$ in bad
$I_{a,q}$'s. Here we use $\sum_{n \leq x} (\frac{n}{\phi(n)})^{1/4}
\leq \sum_{n \leq x} \frac{n}{\phi(n)} \ll x$.

\bigskip

The proof of Theorem \ref{thm3} is almost the same as Theorem 6 of
[\ref{gk}]. So we will give the main points only. The main
modification is the set $\mathcal{V}$ to which we want it to
resemble the set of primes $\leq m^{1/2}$ (the choice in [\ref{gk}])
but each element has size $\approx N$.

\bigskip

Sketch of proof of Theorem \ref{thm3}: Suppose $N \leq m \leq
N^{2-\epsilon}$. Let $I = [N/4, N]$ and $I_k :=
(\frac{m^{1/2}}{2^k}, \frac{m^{1/2}}{2^{k-1}}]$ for $k = 1, 2, 3,
...$ Then $a I_k \subset I$ for $\frac{N 2^k}{4 m^{1/2}} < a \leq
\frac{N 2^k}{2 m^{1/2}}$. Note that $\frac{N}{m^{1/2}} \geq
N^{\epsilon/2}$. By Bertrand's postulate, there is a prime $p_k$
with $\frac{N 2^k}{4 m^{1/2}} < p_k \leq \frac{N 2^k}{2 m^{1/2}}$
and $p_k I_k \subset I$. Now, one of $p_1$, $p_2$, ...,
$p_{[4/\epsilon]+1}$ must be relatively prime to $m$. For otherwise
$p_1 p_2 ... p_{[4/\epsilon]+1}$ divides $m$ which implies
$(\frac{N^{\epsilon/2}}{2})^{[4/\epsilon] + 1} \leq m \leq
N^{2-\epsilon}$. This is impossible for sufficiently large $N$.

Therefore, say for some $1 \leq K \leq [4/\epsilon]+1$, $p_K$ is
relatively prime to $m$. We define $\mathcal{V} := \{ p_K p : p
\hbox{ is a prime in } (\frac{m^{1/2}}{2^K},
\frac{m^{1/2}}{2^{K-1}}], \; (p,m) = 1\}$. Clearly $\mathcal{V}
\subset [N/4, N]$ and $|\mathcal{V}| \gg_\epsilon
\frac{m^{1/2}}{\log m}$ as $m$ is divisible by at most two different
primes in the interval $(\frac{m^{1/2}}{2^K},
\frac{m^{1/2}}{2^{K-1}}]$. Using this $\mathcal{V}$ in the proof of
Theorem 6 of [\ref{gk}], one can get that the set
$$S = \{x y \pmod m : x \in \mathcal{V}, S+1 \leq y \leq S +
\Delta m^{1/2} \sqrt{\frac{m}{\phi(m)}} \log m \}$$ contains $(1 +
O_\epsilon(\Delta^{-1})) m$ residue classes modulo $m$. Now, we want
to add the requirement $(x,y)=1$ by dropping some residue classes.
We are going to exclude those $y$'s that are divisible by $p_K$ or
$p \in (\frac{m^{1/2}}{2^K}, \frac{m^{1/2}}{2^{K-1}}]$. For $p_K$,
we exclude $O(\Delta m^{1/2 - \epsilon} \sqrt{\frac{m}{\phi(m)}}
\log m)$ $y$'s as $p_K \gg N^{\epsilon/2}$. For each $p$, we exclude
$O_\epsilon(\Delta \sqrt{\frac{m}{\phi(m)}} \log m)$ $y$'s. So, in
total, we exclude $O_\epsilon(\Delta m^{1/2}
\sqrt{\frac{m}{\phi(m)}})$ $y$'s. Therefore, we need to exclude at
most $O_\epsilon(\Delta \frac{m}{\log m} \sqrt{\frac{m}{\phi(m)}})$
residue classes from $S$ to ensure $(x,y) = 1$ and we have Theorem
\ref{thm3}.
\section{Erd\"{o}s - Tur\'{a}n inequality: Theorem \ref{thm6}} \label{sec4}

First, we recall a variant of Erd\"{o}s - Tur\'{a}n inequality on
uniform distribution (see, for example, R.C. Baker [\ref{B}, Theorem
2.2]).
\begin{lem} \label{lem1}
Let $L$ and $J$ be positive integers. Suppose that $||x_j|| \geq
\frac{1}{L}$ for $j= 1, 2,... , J$. Then
$$\sum_{l = 1}^{L} \Big| \sum_{j = 1}^{J} e(l x_j) \Big| >
\frac{J}{6}.$$ Here $||x|| = \min_{n \in \mathbb{Z}} |x - n|$, the
distance from $x$ to the nearest integer, and $e(x) = e^{2 \pi i
x}$.
\end{lem}

Proof of Theorem \ref{thm6}: Let $\epsilon > 0$, $N \geq 1$ and $1 +
\epsilon \leq \phi \leq 2$. Suppose $\alpha$ has a rational
approximation $|\alpha - \frac{a}{q}| \leq \frac{1}{q N^{\phi}}$ for
some integers $a$, $N^{2 - \phi + \epsilon} \leq q \leq N^{\phi}$
and $(a,q) = 1$ (Note: The situation when $q < N^{2 - \phi +
\epsilon}$ is trivial as one can simply pick any two distinct primes
$q_1$, $q_2$ in the interval $[N/4, N]$ and successive fractions
with denominator $q_1 q_2$ has spacing $O(\frac{1}{N^2}) =
O(\frac{1}{q N^{2 - \epsilon}})$). We try to find integers $k$ and
distinct prime numbers $q_1, q_2 \in \mathcal{P}$ (the set of prime
numbers in the interval $[N/2, N]$) such that
$$\Big| \frac{a}{q} - \frac{k}{q_1 q_2} \Big| < \frac{1}{q N^{\phi
- \epsilon}} \; \hbox{ or } \; \Big\Vert \frac{q_1 q_2 a}{q}
\Big\Vert < \frac{1}{q N^{\phi - 2 - \epsilon}}.$$ In view of the
above lemma, to prove Theorem \ref{thm6}, it suffices to show
$$\sum_{l = 1}^{L} \Big| \mathop{\sum_{q_1, q_2 \in
\mathcal{P}}}_{q_1 \neq q_2} e \Bigl(\frac{l q_1 q_2 a}{q} \Bigr)
\Big| \leq \frac{1}{6} (|\mathcal{P}|^2 - |\mathcal{P}|)$$ with $L
:= [q N^{\phi - 2 - \epsilon}] + 1$. By triangle inequality, it
suffices to show
\begin{equation} \label{method2}
S_1 + S_2 := \sum_{l = 1}^{L} \Big| \sum_{q_1, q_2 \in \mathcal{P}}
e \Bigl(\frac{l q_1 q_2 a}{q} \Bigr) \Big| + \sum_{l = 1}^{L} \Big|
\sum_{q_1 \in \mathcal{P}} e \Bigl(\frac{l q_1^2 a}{q} \Bigr) \Big|
\leq \frac{1}{7} |\mathcal{P}|^2
\end{equation}
for $N$ sufficiently large. Note that $1 \leq L \leq q$ as $1 +
\epsilon \leq \phi \leq 2$. By Cauchy-Schwarz inequality and
orthogonality of $e(x)$,
\begin{equation*}
\begin{split}
S_2^2 \leq& \Bigl(\sum_{l = 1}^{L} 1 \Bigr) \Big(\sum_{l = 1}^{L}
\Big|\sum_{q_1 \in \mathcal{P}} e \Bigl(\frac{l q_1^2 a}{q} \Bigr)
\Big|^2 \Bigl) \leq L \Big(\sum_{l = 1}^{q} \Big|\sum_{q_1 \in
\mathcal{P}} e \Bigl(\frac{l q_1^2 a}{q} \Bigr)
\Big|^2 \Bigl) \\
=& L q \mathop{\sum_{q_1 \in \mathcal{P}} \sum_{q_2 \in
\mathcal{P}}}_{q_2^2 \equiv q_1^2 \pmod q} 1 \ll \left\{
\begin{tabular}{ll} $L q d(q) N$ & if
$q > N$, \\
$L d(q) N^2$ & if $q \leq N$
\end{tabular} \right.
\end{split}
\end{equation*}
because $(q_1 ,q) = 1 = (q_2,q)$, and $n^2 \equiv b \pmod q$ has
$O(2^{\omega(q)}) = O(d(q))$ solutions for $n$ when $(b,q)=1$. Here
$\omega(q)$ denotes the number of distinct prime divisors of $q$.
Therefore as $L = [q N^{\phi - 2 - \epsilon}] + 1$ and $q \leq
N^\phi$,
\begin{equation*}
\begin{split}
S_2 \ll& \left\{ \begin{tabular}{ll} $q N^{\phi/2 - 1/2 -
\epsilon/2} d(q)^{1/2}$ & if $q > N$, \\
$q^{1/2} N^{\phi/2 - \epsilon/2} d(q)^{1/2}$ & if $q \leq N$,
\end{tabular} \right. \\
\ll& d(q)^{1/2} \max(N^{3\phi/2 - 1/2 - \epsilon/2}, N^{\phi/2 + 1/2
- \epsilon/2}).
\end{split}
\end{equation*}
By Chebychev's estimate or the prime number theorem, and $d(q)
\ll_\epsilon q^{\epsilon/4} \leq N^{\epsilon / 2}$, one can check
that $S_2 \leq \frac{1}{14} |\mathcal{P}|^2$ when $\phi \leq 5/3$
and $N$ is sufficiently large depending on $\epsilon$.

It remains to deal with $S_1$. Our approach is inspired by
Vinogradov's work [\ref{V}]. Write
$$d_r = \mathop{\sum_{l \leq L, q_1 \in \mathcal{P}}}_{l q_1 = r}
1.$$ Since $l \leq L \leq q N^{\phi - 2 - \epsilon} + 1 \leq
N^{2\phi - 2 - \epsilon} + 1 \leq N^{2 - \epsilon} + 1 \ll N^2$, we
need to consider $r \leq LN \ll N^3$. But then $r$ is divisible by
at most three $q_1 \in \mathcal{P}$. Thus $d_r \leq 3$. Then
\begin{equation} \label{s1}
S_1 \leq \sum_{r = 1}^{LN} d_r \Big| \sum_{q_2 \in \mathcal{P}} e
\Bigl(\frac{r q_2 a}{q}\Bigr) \Big|
\end{equation}
By Cauchy-Schwarz inequality and orthogonality of $e(x)$,
\begin{equation*}
\begin{split}
S_1^2 \leq& \Bigl( \sum_{r = 1}^{LN} d_r^2 \Bigr)
\Bigl(\sum_{r=1}^{LN} \Big| \sum_{q_2 \in \mathcal{P}} e
\Bigl(\frac{r q_2 a}{q}\Bigr) \Big|^2 \Bigr) \leq 9LN \sum_{r =
1}^{q([LN/q]+1)} \Big| \sum_{q_2 \in \mathcal{P}} e
\Bigl(\frac{r q_2 a}{q}\Bigr) \Big|^2 \\
=& 9LN \Bigl([\frac{LN}{q}] + 1 \Bigr) q \mathop{\sum_{q_1 \in
\mathcal{P}} \sum_{q_2 \in \mathcal{P}}}_{q_2 \equiv q_1 \pmod q} 1 \\
\leq& 9LN \Bigl([\frac{LN}{q}] + 1 \Bigr) q \mathop{\sum_{N/2 \leq
q_1 \leq N} \sum_{N/2 \leq q_2 \leq N}}_{q_2 \equiv q_1 \pmod q} 1
\ll \left\{ \begin{tabular}{ll} $(L N)^2 N$ & if
$q > N$, \\
$(L N)^2 \frac{N^2}{q}$ & if $q \leq N$
\end{tabular} \right.
\end{split}
\end{equation*}
as $LN \geq q N^{\phi - 1 - \epsilon} \geq q$. Therefore as $L = [q
N^{\phi - 2 - \epsilon}] + 1$ and $q \leq N^\phi$,
$$S_1 \ll \left\{ \begin{tabular}{ll} $q N^{\phi - 1/2 -
\epsilon}$ & if $q > N$ \\
$q^{1/2} N^{\phi - \epsilon}$ & if $q \leq N$
\end{tabular} \right. \ll \max (N^{2\phi - 1/2 - \epsilon}, N^{\phi
+ 1/2 - \epsilon}).$$

By Chebychev's estimate or the prime number theorem, one can check
that $S_1 \leq \frac{1}{14} |\mathcal{P}|^2$ when $\phi \leq 5/4$
and $N$ is sufficiently large depending on $\epsilon$. Consequently,
we have (\ref{method2}) as long as $\phi \leq 5/4$ and $N$
sufficiently large. We set $\phi = 5/4$.

By the contrapositive of Lemma \ref{lem1}, there exist distinct
prime numbers $q_1, q_2 \in \mathcal{P}$ such that $||\frac{q_1 q_2
a}{q}|| < \frac{1}{q N^{5/4 - 2 - \epsilon}}$ for sufficiently large
$N$. In other words $|\frac{a}{q} - \frac{k}{q_1 q_2}| \ll_\epsilon
\frac{1}{q N^{5/4 - \epsilon}}$ for some integer $k$. Hence
$$\Big| \alpha - \frac{k}{q_1 q_2} \Big| \leq \Big|\alpha -
\frac{a}{q} \Big| + \Big|\frac{a}{q} - \frac{k}{q_1 q_2} \Big|
\ll_\epsilon \frac{1}{q N^{5/4 - \epsilon}}$$ which gives Theorem
\ref{thm6} as $(q_1, q_2) = 1$ and we can write $k = a_1 q_2 + a_2
q_1$ for some integers $a_1$, $a_2$.
\section{Theorem \ref{thm7}} \label{sec5}
The proof of Theorem \ref{thm7} is almost the same as Theorem
\ref{thm6}. We shall be content to indicate the necessary
modifications.

\bigskip

Proof of Theorem \ref{thm7}: Without loss of generality, we can
assume $q > N$ as indicated in the note of Conjecture \ref{conj1}.
The starting point is almost the same as that of Theorem \ref{thm6}.
The difference is that we want $q_1 \in \mathcal{P}$ but any integer
$q_2 \in [N/2, N]$ with $(q_1, q_2) = 1$. Equivalently we need $q_2
\neq q_1$ since $q_1$ is a prime in $[N/2, N]$. Instead of
(\ref{method2}), it suffices to show
\begin{equation} \label{method3}
S_1 + S_2 := \sum_{l = 1}^{L} \Big| \sum_{q_1 \in \mathcal{P}}
\sum_{N/2 \leq q_2 \leq N} e \Bigl(\frac{l q_1 q_2 a}{q} \Bigr)
\Big| + \sum_{l = 1}^{L} \Big| \sum_{q_1 \in \mathcal{P}} e
\Bigl(\frac{l q_1^2 a}{q} \Bigr) \Big| \leq \frac{1}{7}
|\mathcal{P}|^2
\end{equation}
Now $S_2$ is treated the same way as in Theorem \ref{thm6}. So $S_2
\leq \frac{1}{14} |\mathcal{P}|^2$ when $\phi \leq 5/3$ and $N$ is
sufficiently large depending on $\epsilon$.

As for $S_1$, instead of (\ref{s1}), we have
$$S_1 \leq \sum_{r = 1}^{q([LN/q]+1)} d_r \Big| \sum_{N/2 \leq q_2 \leq N} e
\Bigl(\frac{r q_2 a}{q}\Bigr) \Big| \ll \frac{LN}{q} \sum_{r =
1}^{q} \Big| \sum_{N/2 \leq q_2 \leq N} e \Bigl(\frac{r
q_2}{q}\Bigr) \Big|$$ as $(a,q) = 1$. Now summing according to the
greatest common divisor $d = (r,q)$,
\begin{equation*}
\begin{split}
S_1 \ll& \frac{LN}{q} \sum_{d | q} \mathop{\sum_{r' \leq q/d}}_{(r',
q/d) = 1} \Big| \sum_{N/2 \leq q_2 \leq N} e \Bigl(\frac{r'
q_2}{q/d}\Bigr) \Big| \leq \frac{LN}{q} \sum_{d | q} \sum_{r' \leq
q/d} \min \Bigl(N,
\frac{1}{||\frac{r'}{q/d}||} \Bigr) \\
\ll& \frac{LN}{q} \mathop{\sum_{d | q}}_{d < N} \Bigl( \sum_{r' \leq
q/N} N + \sum_{q/N < r' \leq q/d} \frac{q/d}{r'} \Bigr) +
\frac{LN}{q} \mathop{\sum_{d | q}}_{d \geq N} \sum_{r' \leq q/d} N
\ll L N d(q) \log q
\end{split}
\end{equation*}
Recall $L = [q N^{\phi - 2 - \epsilon}] + 1$ and $q \leq N^\phi$,
thus $S_1 \leq \frac{1}{14} |\mathcal{P}|^2$ as long as $\phi \leq
3/2 < 5/3$. This proves Theorem \ref{thm7}.
\section{Character sum method: Theorem \ref{thm5}} \label{sec3}
First let us prove a simple lemma which is needed at the end of the
proof of Theorem \ref{thm5}.
\begin{lem} \label{charlem}
Let $q$ be a positive integer and $B \geq 1$ be any real number.
Then the number of integers between $1$ and $B$ relatively prime to
$q$ is $B \frac{\phi(q)}{q} + O(d(q))$.
\end{lem}

Proof: Using properties of M\"{o}bius function $\mu(n)$, the number
of integers between $1$ and $B$ relatively prime to $q$ is
\begin{equation*}
\begin{split}
\mathop{\sum_{1 \leq n \leq B}}_{(n,q) = 1} 1 =& \sum_{1 \leq n \leq
B} \mathop{\sum_{d | n}}_{d | q} \mu(d) = \mathop{\sum_{1 \leq d
\leq B}}_{d | q} \mu(d) \mathop{\sum_{1 \leq n \leq B}}_{d | n} 1 =
\mathop{\sum_{1 \leq d \leq B}}_{d | q} \mu(d) \Bigl(\frac{B}{d}
+ O(1) \Bigr) \\
=& B \sum_{d | q} \frac{\mu(d)}{d} + O \Bigl(B \mathop{\sum_{d >
B}}_{d | q} \frac{\mu(d)}{d} \Bigr) + O(d(q)) = B \frac{\phi(q)}{q}
+ O(d(q)).
\end{split}
\end{equation*}

Proof of Theorem \ref{thm5}: Let $\epsilon > 0$, $N \geq 1$ and $1 +
\epsilon \leq \phi \leq 2 - 2\epsilon$. Without loss of generality,
we may assume $N$ is sufficient large. Suppose $\alpha$ has a
rational approximation $|\alpha - \frac{a}{q}| \leq \frac{1}{q
N^{\phi}}$ for some integers $a$, $N^{2 - \phi - \epsilon} < q \leq
N^{\phi}$ and $(a,q) = 1$ (the case $q \leq N^{2 - \phi - \epsilon}$
is trivial as indicated in the proof of Theorem \ref{thm6}). We will
try to find integer $k$ and distinct prime numbers $N/2 \leq q_1,
q_2 \leq N$ such that $|\frac{a}{q} - \frac{k}{q_1 q_2}|$ is small.
This is equivalent to $|a q_1 q_2 - k q|$ being small which leads us
to consider
\begin{equation} \label{cong}
a q_1 q_2 \equiv b \pmod q
\end{equation}
with distinct $q_1, q_2 \in \mathcal{P}$, the set of primes in the
interval $[N/2, N]$ that are relatively prime to $q$, and $b \in
\mathcal{B}$, the set of integers in the interval $[1,B]$ with $B
\leq q$ to be chosen later.

Let $\chi$ denote a typical Dirichlet character modulo $q$. By
orthogonality of characters, the number of solution to (\ref{cong})
is
\begin{equation} \label{number}
\# = \frac{1}{\phi(q)} \sum_{\chi} \mathop{\sum_{q_1, q_2 \in
\mathcal{P}}}_{q_1 \neq q_2} \sum_{b \in \mathcal{B}} \chi(a q_1
q_2) \overline{\chi}(b)
\end{equation}
where the sum $\sum_\chi$ is over all Dirichlet character modulo
$q$, and $\overline{z}$ denotes complex conjugate of $z$. We want
$\# > 0$. We separate the contribution from the principal character
in (\ref{number}) and get
$$\# = \frac{1}{\phi(q)} |\mathcal{P}|(|\mathcal{P}| - 1) |\mathcal{B}_q|
+ \frac{1}{\phi(q)} \sum_{\chi \neq \chi_0} \mathop{\sum_{q_1, q_2
\in \mathcal{P}}}_{q_1 \neq q_2} \sum_{b \in \mathcal{B}} \chi(a q_1
q_2) \overline{\chi}(b) := \#_1 + \#_2$$ where $\mathcal{B}_q$
denotes the set of numbers in $\mathcal{B}$ that are relatively
prime to $q$. Now
\begin{equation*}
\begin{split}
\#_2 =& \frac{1}{\phi(q)} \sum_{\chi \neq \chi_0} \sum_{q_1, q_2 \in
\mathcal{P}} \sum_{b \in \mathcal{B}} \chi(a q_1 q_2)
\overline{\chi}(b) - \frac{1}{\phi(q)} \sum_{\chi \neq \chi_0}
\sum_{q_1 \in \mathcal{P}} \sum_{b \in \mathcal{B}} \chi(a q_1^2)
\overline{\chi}(b) \\
|\#_2| \leq& \frac{1}{\phi(q)} \sum_{\chi \neq \chi_0} \Big|
\sum_{q_1 \in \mathcal{P}} \chi(q_1) \Big|^2 \Big|\sum_{b \in
\mathcal{B}} \chi(b) \Big| + \frac{1}{\phi(q)} \sum_{\chi \neq
\chi_0} \Big| \sum_{q_1 \in \mathcal{P}} \chi^2(q_1) \Big|
\Big|\sum_{b \in \mathcal{B}} \chi(b) \Big|
\\
:=& S_1 + S_2.
\end{split}
\end{equation*}
Using trivial estimate on $\sum_{q_1}$ and Cauchy-Schwarz
inequality, we have
$$S_2 \leq \frac{|\mathcal{P}|}{\phi(q)} \Bigl(\sum_{\chi} 1 \Bigr)^{1/2}
\Bigl(\sum_{\chi} \Big| \sum_{b \in \mathcal{B}} \chi(b) \Big|^2
\Bigr)^{1/2} \leq |\mathcal{P}| |\mathcal{B}|^{1/2} \leq
|\mathcal{B}|^{1/2} N$$ by orthogonality of characters (note that $B
\leq q$).

Now, recall a well-known consequence of the Generalized Lindel\"{o}f
Hypothesis:
\begin{equation} \label{GLH}
\sum_{n \leq N} \chi(n) \ll_\epsilon N^{1/2} q^{\epsilon}
\end{equation}
for non-principal character $\chi \pmod q$. Using this and the
orthogonality of characters, we have
\begin{equation*}
\begin{split}
S_1 \ll_\epsilon& |\mathcal{B}|^{1/2} q^\epsilon \frac{1}{\phi(q)}
\sum_{\chi} \Big| \sum_{q_1 \in \mathcal{P}} \chi(q_1) \Big|^2 \leq
|\mathcal{B}|^{1/2} q^\epsilon \mathop{\sum_{q_1, q_1' \in
\mathcal{P}}}_{q_1 \equiv q_1' \pmod q} 1 \\
\leq& |\mathcal{B}|^{1/2} q^\epsilon \mathop{\sum_{N/2 \leq q_1,
q_1' \leq N}}_{q_1 \equiv q_1' \pmod q} 1 \leq |\mathcal{B}|^{1/2}
q^\epsilon N (\frac{N}{q} + 1).
\end{split}
\end{equation*}

Thus
$$|\#_2| \ll_\epsilon |\mathcal{B}|^{1/2}
q^\epsilon N (\frac{N}{q} + 1) \ll \bigl\{ \begin{tabular}{ll}
$|\mathcal{B}|^{1/2} q^\epsilon N$ & if $q > N$, \\
$|\mathcal{B}|^{1/2} q^\epsilon \frac{N^2}{q}$ & if $q \leq N$;
\end{tabular}$$ while
$$\#_1 \geq \frac{|\mathcal{P}|^2 |\mathcal{B}_q|}{2 \phi(q)}$$
for $N$ sufficiently large. Therefore, replacing $\epsilon$ with
$\epsilon/4$, we have $\# > 0$ if
$$\frac{|\mathcal{P}|^2 |\mathcal{B}_q|}{\phi(q)} \gg_\epsilon
|\mathcal{B}|^{1/2} q^{\epsilon/4} N \; \hbox{ and } \;
\frac{|\mathcal{P}|^2 |\mathcal{B}_q|}{\phi(q)} \gg_\epsilon
|\mathcal{B}|^{1/2} \frac{N^2}{q^{1 - \epsilon/4}}.$$ By Chebyshev's
estimate or the prime number theorem, $|\mathcal{P}| \geq \frac{N}{2
\log N} - 2 > \frac{N}{3 \log N}$ as $q \leq N^2$ can be divisible
by at most two $q_1 \in \mathcal{P}$. Consequently by Lemma
\ref{charlem}, (\ref{cong}) has some solutions $q_1, q_2 \in
\mathcal{P}$ with $q_1 \neq q_2$, and $1 \leq b \leq B =
\max(\frac{q^{2 + 0.6 \epsilon}}{N^2}, q^{0.6 \epsilon})$. In other
words, we can find distinct primes $q_1 ,q_2 \in \mathcal{P}$ and
integer $k$ such that
\begin{equation} \label{3bound1}
\Big| \frac{a}{q} - \frac{k}{q_1 q_2} \Big| \ll_\epsilon \max \Bigl(
\frac{q^{1 + 0.6 \epsilon}}{N^4}, \frac{1}{q^{1 - 0.6 \epsilon} N^2}
\Bigr).
\end{equation}
The right hand side of (\ref{3bound1}) is $\ll_\epsilon \frac{1}{q
N^{\phi - \epsilon}}$ if $q^{2 + 0.6\epsilon} \leq N^{4 - \phi +
\epsilon}$. This is true when $\phi \leq 4/3$ as $q \leq N^\phi$
(Note: $0.6 \times 4/3 < 1$). This proves Theorem \ref{thm5}.

Note: One can get a weaker unconditional result with $\phi \leq 6/5$
using P\'{o}lya-Vinogradov inequality instead of (\ref{GLH}).

\bigskip

{\bf Acknowledgement} The author would like to thank Angel Kumchev
for some helpful comments on the almost all results, Sidney Graham
for some stimulating discussions, and George Grossman for
suggestions on exposition. The author also thanks the referee for
some helpful corrections.

\bigskip
Department of Mathematical Sciences \\
University of Memphis \\
Memphis, TN 38152 \\
U.S.A. \\
tszchan@memphis.edu
\end{document}